\documentclass[12pt]{article}
\usepackage{amsfonts,amssymb}
\usepackage{amsmath}
\usepackage{latexsym}  
\usepackage{authblk}

\title{Density results and 
trace operator in weighted Sobolev Spaces defined on the half line equipped with power weights }
\author{{\large{Rados{\l}aw Kaczmarek}}\vspace{-0.5cm}\\
\textit{Faculty of Mathematics and Computer Science, Adam Mickiewicz University in Pozna\'{n}, Uniwersytetu Pozna\'{n}skiego 4, 61-614 Pozna\'{n}, Poland; radekk@amu.edu.pl, ORCID: 0000-0002-3536-9102} 
\vspace{0.5cm}\\
{\large{Agnieszka Ka{\l}amajska}}\\
\textit{Faculty of Mathematics, Informatics and Mechanics, University of Warsaw, Banacha 2, 02-097 Warsaw, Poland; A.Kalamajska@mimuw.edu.pl, ORCID: 0000-0001-5674-8059}}

\usepackage{color}

\newtheorem{theo}{\bf Theorem}[section]

\newtheorem{rem}{\bf Remark}[section]

\newtheorem{lemma}{\bf Lemma}[section]

\makeatletter \@addtoreset{equation}{section}

\makeatother

\newcommand{\beq}{\begin{equation}}
\newcommand{\eeq}{\end{equation}}

\newcommand{\BR}{\mathbb{R}}
\newcommand{\BN}{\mathbb{N}}

\newcommand{\ve}{\varepsilon}

\newcommand{\rp}{\mathbb{R}_+}

\newcommand{\ypb}{Y^{1,p}_{t^\beta}}

\begin{document}

\maketitle

\begin{abstract}
We study properties of $W_0^{1,p}(\rp,t^\beta)$ - the completion of $C_0^\infty(\BR_+)$ in the power-weighted Sobolev spaces $W^{1,p}(\BR_+,t^\beta)$,  where $\beta\in\BR$.
Among other results, we obtain the analytic characterization of  $W_0^{1,p}(\rp,t^\beta)$ for all $\beta\in \mathbb{R}$. Our analysis is based on the precise study the two trace operators: $Tr^{0}(u):= \lim_{t\to 0} u(t)$ and 
$Tr^{\infty}(u):= \lim_{t\to \infty} u(t)$, which leads to the  analysis of the  asymptotic behavior of functions from $W_0^{1,p}(\rp,t^\beta)$ near zero or infinity. The obtained statements can contribute to the proper formulation of Boundary Value Problems in ODE's, or PDE's with the radial symmetries.
We  can also apply our results to some questions in the complex  interpolation theory, raised by M. Cwikel and A. Einav in 2019, which we discuss within the particular case of Sobolev spaces  $W^{1,p}(\BR_+,t^\beta)$. 
\end{abstract}

\begin{flushleft}
\end{flushleft}

\section{Introduction}
Results about densities of smooth functions are one of central problems in the theory of Sobolev spaces, see e.g. \cite{af,ma}. Classically,
having the domain $U\subseteq\mathbb{R}^n$, one defines $W_0^{m,p}(U)$ - the subspace of Sobolev space $W^{m,p}(U)$, as the completion of the set of  smooth compactly supported  functions  $C_0^\infty(U)$, in the norm of  $W^{m,p}(U)$. Similarly, having given the weights function $\omega$ defined on $U$,  one defines $W_0^{m,p}(U,\omega)$ as  the completion of $C_0^\infty (U)$ in the weighted Sobolev space $W^{m,p}(U,\omega)$. Then the following problem arises: {\it How can we recognize in the easy way, if the given function from the space $W^{m,p}(U,\omega)$ belongs to its subspace $W_0^{m,p}(U,\omega)$?}

Here we are interested in spaces  $W^{1,p}(\rp,t^\beta)$ and $W_0^{1,p}(\rp,t^\beta)$, where $\beta\in\mathbb{R}$ is given number and we ask about the analytic characterization of $W_0^{1,p}(\rp,t^\beta)$. Although the above question looks very simple at first glance, we could not find results about such a characterization in the literature.  
On the other hand, they might be present in the hidden way, as it is in the classical formulation of Hardy inequality, see Remark \ref{remhardy}. 

\smallskip

In the non-weighted, classical approach on bounded and sufficiently regular domains, characterization is known. Namely, having given bounded domain with the sufficiently regular boundary, say Lipschitz, one has well defined the trace operator $Tr : W^{1,p}(U)\rightarrow L^1(\partial U)$ (see e.g. \cite{af,kjf,ma}) for any $p\in [1,+\infty)$. Then the following characterization can be prescribed:
$$
W_0^{1,p}(U)= \{ u\in W^{1,p}(U): Tr \;u|_{\partial U}\equiv 0\ {a.e.}\ {\rm on}\  \partial U\}.
$$ 
Our analysis is based on a similar approach, which however requires
to have introduced and analyzed
the {\it  trace operator} defined at the endpoint of the interval $(0,\infty)$. It appears, that in our setting, in all cases except at $\beta =p-1$, such trace operator can be defined at one endpoint $p=0$ or $p=\infty$ only. For example, as shown in Theorems: \ref{pierwsze}, \ref{trzecie} and \ref{drugie}, for any $p>1$ and $u\in W^{1,p}(\rp, t^\beta)$, we have:
\begin{itemize}
\item $Tr^0(u):=\lim_{t\to 0}u(t)$ is well prescribed for $\beta<p-1$,
\item $Tr^\infty (u):=\lim_{t\to \infty}u(t)$ is well prescribed for $\beta>p-1$,
\item when $\beta = p-1$, then $Tr^0(u)$ cannot be prescribed in general on $W^{1,p}(\rp,t^{p-1})$, while $Tr^\infty(u)$ annihilates all that space.
\end{itemize}

Moreover,  having trace operator precisely defined,  when $u\in W_0^{1,p}(\rp,t^\beta)$, then one has
 \begin{eqnarray*}
 \lim_{t\to 0} u(t)= 0 \Longleftrightarrow \lim_{t\to 0} t^{(\beta-(p-1))/p} u(t),\ \ \hbox{\rm when}\ \beta <p-1,\\
  \lim_{t\to \infty} u(t)= 0 \Longleftrightarrow \lim_{t\to \infty} t^{(\beta-(p-1))/p} u(t),\ \ \hbox{\rm when}\ \beta >p-1,
  \end{eqnarray*}
which contributes to the study of the {\it asymptotic behavior} of functions from $W_0^{1,p}(\rp,t^\beta)$, see Theorems \ref{pierwsze} and \ref{trzecie}. The above facts can be applied to  proper {\it formulation of Boundary Value Problems} in ODE's and radial solutions to PDE's with radial symmetries, as we mention in Remark \ref{applbvp}.
 The case of $\beta = p-1$ is the borderline one between that of $\beta <p-1$ and of $\beta>p-1$. It appears that in that case  $W_0^{1,p}(\rp, t^{p-1})=W^{1,p}(\rp, t^{p-1})$, which is shown in Theorem \ref{drugie}. As just mentioned, the trace operators do not apply reasonably to that space. 

\smallskip
New information about asymptotic behavior can be also applied to {\it Hardy type inequalities}. Namely, 
the  set $W=W(\alpha,\beta)$ such that $C_0^\infty(\rp)\subseteq W\subseteq W^{1,p}(\rp,t^\beta)$ and $W$ is admitted to the inequality:
$$\int_{\rp} |t^\alpha |u(t)||^pt^\beta dt \le C_1 \int_{\rp} |u(t)|^pt^\beta dt + C_2\int_{\rp} |u^{'}(t)|^pt^\beta dt, $$
has been precisely described in the paper \cite{akikppcejm} in the purely analytic way. Now we can say that this set is precisely $W_0^{1,p}(\rp,t^\beta)$ in the case when $\beta < |\alpha |p -1$, $\alpha \in (-1,0)$, which is not visible without our analysis. See Remark \ref{appl-interpo2} for more precise information.

\smallskip

Other problems of our interest, are {\it density results}. Having traces well defined, we can now say that:
\begin{itemize}
\item
$W_0^{1,p}(\rp, t^\beta) = \{ u\in W^{1,p}(\rp, t^\beta): \lim_{t\to 0}u(t)=0\}$ when $\beta<p-1$, 
\item
$W_0^{1,p}(\rp, t^\beta) = \{ u\in W^{1,p}(\rp, t^\beta): \lim_{t\to \infty}u(t)=0\}$ when $\beta>p-1$,
\item
$W_0^{1,p}(\rp, t^\beta) = W^{1,p}(\rp, t^\beta)$, when $\beta = p-1$,
\end{itemize} 
see Theorems: \ref{pierwsze}, \ref{trzecie}, \ref{drugie}.

As next example application of our results, we focus on the {\it application to  complex interpolation theory}. Namely, our density results allow us to conclude that
$$
 [W^{1,p}(\rp,t^{\beta_0},\mathbb{C}), W^{1,p}(\rp,\omega_1, \mathbb{C})]_\theta 
= W^{1,p}(\rp ,t^{\beta_\theta},  \mathbb{C}),
$$
whenever, $\beta_\theta:= (1-\theta)\beta_0 + \theta\beta_1\in (-\infty,-1]\cup (p-1,\infty)$ and $[X,Y]_\theta$ is complex interpolation pair between Banach spaces $X,Y$, with $\theta\in (0,1)$, see Theorem \ref{interpo-konkluzja-my}. Such result seems new and contributes to the problem posed in 1982 by L\"ofstr\"om in  \cite{Lo} and to the discussion undertaken  recently by M. Cwikel and A. Einav in \cite{ce} about interpolation spaces between weighted Sobolev spaces. 
In particular, we built some examples partially solving  questions posed in \cite{ce}, when restricting the analysis to the particular cases of weighted Sobolev spaces $W^{1,p}(\BR_+,t^\beta)$.  See Remark \ref{appl-interpo}  for the precise  formulations and arguments.

\smallskip

{\it The methods} we use are rather elementary. In most situations they are based on the appriori estimates, first order integral Taylor's formulae, and on Hardy inequality (Theorem \ref{hardyclas}). Similar tools were used e.g. in \cite{akikpp2009} on page 9, as well as in \cite{kufner}, in the proof of Lemma 5.9. More abstract tools were presented in the proof of Theorem \ref{drugie}, showing that $W^{1,p}(\BR_+,t^\beta)= W^{1,p}_0(\BR_+,t^\beta)$, where we used Mazur's Lemma.

\smallskip

As all of them: density results, asymptotic behavior near endpoints as well as complex interpolation theory are  important tools in functional analysis and PDE's/ODE's, we hope to contribute in this way for their understanding.

\section{Notation and Preliminaries}

{\bf Basic notation.}
We will deal with functions defined on  $\rp$  and use the standard notations $C^\infty(\BR_+)$, $C^\infty_0(\BR_+)$, ${Lip}(\BR_+)$, $L^p(\rp,\omega)$ for smooth functions, smooth compactly supported functions, lipschitz functions and for functions in $L^p$ equipped with $\omega$ weight, respectively. 
 If $A$ is a subset in some Banach space $X=(X,\|.\|_X)$, then by $\overline{A}^{\,(\|.\|_X)}$ we denote the completion of $A$ in $X$ in the norm $\|.\|_X$. 

When estimating expressions, notation $A\preceq B$ will be used when the estimate $A\leq CB$ holds with some universal constant $C$, whose precise value  is not important for final conclusion, while $A\sim B$ will mean that 
$A\preceq B$ and $B\preceq A$ (denoted also as $A \succsim B$).

\bigskip
\noindent
{\bf The special Sobolev spaces.}
Let $1\leq p<\infty$ and $\beta\in\BR$. We will deal with the following one and two weighted Sobolev spaces $W^{1,p}(\BR_+,t^\beta )$ and $Y^{1,p}_{t^\beta}$:
\begin{eqnarray}\label{sobolewy1} 
W^{1,p}(\BR_+,t^\beta ) &:=&
\{u\in W^{1,1}_{loc}(\BR_+)\colon \|u\|_{W^{1,p}(\BR_+,t^\beta)}<\infty \} \ {\rm where}\\
\|u\|_{W^{1,p}(\BR_+,t^\beta)} &:=&\left(\int_{\BR_+} |u(t)|^p t^\beta dt\right)^{\frac{1}{p}}+\left(\int_{\BR_+} |u^\prime (t)|^p t^\beta dt\right)^{\frac{1}{p}}, \nonumber\\
\label{sobolewy2} 
Y^{1,p}_{t^\beta} &:=& \{u\in W^{1,1}_{loc}(\BR_+)\colon \|u\|_{Y^{1,p}_{t^\beta}}<\infty \}  \ {\rm where}\\
\|u\|_{Y^{1,p}_{t^\beta}} &:=&\left(\int_{\BR_+} |u(t)|^p t^\beta (1+\frac{1}{t^p}) dt\right)^{\frac{1}{p}}+\left(\int_{\BR_+} |u^\prime (t)|^p t^\beta dt\right)^{\frac{1}{p}}. \nonumber
\end{eqnarray}
Those spaces are complete and  they form  Banach spaces, see e.g. \cite{kuf-opic}, Theorem 1.11.  

By $W^{1,p}_0(\rp , t^\beta)$, $Y^{1,p}_{t^\beta ,0}$, we denote the completion of $C_0^\infty (\rp )$ in $W^{1,p}(\rp , t^\beta)$ and $Y^{1,p}_{t^\beta }$, respectively. We will also use notation ${Lip_{c}}(\BR_+)$, $W^{1,p}_{c}(\rp,t^\beta)$, 
 $Y^{1,p}_{t^\beta ,c }$ for the sets of functions in a given function space, which have compact support in $\rp$.

Let us recall that $W^{1,1}_{loc}(\rp)\subseteq C(\rp)$, see e.g. \cite{ma}, therefore all the functions which belong to weighted Sobolev spaces as in \eqref{sobolewy1} and \eqref{sobolewy2}, are continuous on $\rp$. 

\bigskip
\noindent
{\bf Hardy inequality.}
We recall the classical Hardy inequality (see e.g. \cite{hlp}, Theorem 330 or \cite{kmp}), which  will be useful for further considerations.

\begin{theo}[Hardy inequality]
\label{hardyclas}
Let $1<p<\infty$, $\beta\neq p-1$. Suppose that $u=u(t)$ is an
absolutely continuous function in $(0,\infty)$ such that
$\int_{0}^\infty |u'(t)|^p\,t^{\beta}dt <\infty,$ and let
\begin{eqnarray*}
u(0)&:=& \lim_{t\to 0} u(t) =0\ \ {\rm for}\ \ \beta <p-1,\\
u(\infty)&:=& \lim_{t\to\infty} u(t)=0\ \ {\rm for}\ \ \beta
>p-1.
\end{eqnarray*}
Then the following inequality holds with sharp constant:
\begin{equation}
\label{clasha}
\int_{0}^\infty |u(t)|^pt^{\beta -p}dt \le \left( \frac{p}{|\beta -p +1| }\right)^p\int_0^\infty
|u'(t)|^pt^{\beta}dt . 
\end{equation}
\end{theo}

\section{First density result}
We start with the following preliminary density result, which deals with Sobolev spaces  as in \eqref{sobolewy1} and \eqref{sobolewy2}. 
 
\begin{lemma} 
\label{gestosci}
Let $1<p<\infty$ and $\beta\in\BR$. 
Then
\begin{eqnarray}\label{jed-a}
W^{1,p}_0(\rp , t^{\beta})=\overline{Lip_c (\rp)}^{\,(\|.\|_{W^{1,p}(\rp , t^{\beta})})}= \overline{W^{1,p}_{c}(\rp , t^{\beta})}^{\,(\|.\|_{W^{1,p}(\rp , t^{\beta})})},
\end{eqnarray}
\begin{eqnarray}\label{jed-b}
Y^{1,p}_{t^{\beta}, 0}= \overline{Lip_c(\rp)}^{\,(\|.\|_{Y^{1,p}_{t^{\beta} }})} = \overline{Y^{1,p}_{t^{\beta},c}}^{\,(\|.\|_{Y^{1,p}_{t^\beta}})} 
\end{eqnarray}
and
\begin{eqnarray}
 Y^{1,p}_{{t^\beta}, 0}= \ypb \subseteq W^{1,p}_0(\rp , t^{\beta}).\label{jed-c}
\end{eqnarray}
\end{lemma}

\begin{rem}\rm 
We will show in the next section that the embedding $\ypb \subseteq W^{1,p}_0(\rp , t^{\beta})$ in \eqref{jed-c} can be either strict or non-strict, depending on the exponent  $\beta$. For example, if $\beta\neq p-1$, then $\ypb=W^{1,p}_0(\rp , t^{\beta})$, while for $\beta=p-1$, we have $\ypb\varsubsetneq W^{1,p}_0(\rp , t^{\beta})$. See Theorems \ref{pierwsze}, \ref{trzecie} and \ref{drugie}.
\end{rem}

\noindent
{\bf Proof.}
{\bf (\ref{jed-a}):} Clearly,
$$
C_0^\infty (\rp) \subseteq Lip_c (\rp ) \subseteq W^{1,p}_{c}(\rp ,t^\beta).
$$
It suffices to show that $W^{1,p}_{c}(\rp ,t^\beta)\subseteq W^{1,p}_{0}(\rp ,t^\beta)$. For that, take $u\in W^{1,p}_{c}(\rp ,t^\beta)$ with the support $[a,b]\subseteq \rp$.
As on compactly supported sets $t^\beta \sim 1$, 
therefore $u\in W^{1,p}(\rp)$, and  is compactly supported. By the standard convolution arguments, the convolutions  
$u_\varepsilon (x):= \phi_\varepsilon *u$ with the classical mollifier functions $\phi_\varepsilon (x)=\varepsilon^{-1}\phi (x/\varepsilon)$  where 
$\phi\in C_0^\infty (\mathbb{R})$, $0\le \phi\le 1$, ${\rm supp\,}u\subseteq [-1,1]$ and $\int_{\mathbb{R}} \phi dx =1$, converge to $u$ in $W^{1,p}(\rp)$. Their supports are the subsets of $J:=[a-\ve_0, b+\ve_0]$ for $\ve_0 := a/2$ and $\varepsilon\in(0,a/2)$. Again, as $t^\beta\sim 1$  on $J$, therefore  $u_\varepsilon$ converge to $u$ also in $W^{1,p}(\rp ,t^\beta)$, as $\varepsilon \to 0$.
This shows that $u\in W^{1,p}_{0}(\rp ,t^\beta)$. 

\smallskip
\noindent
{\bf (\ref{jed-b}):}  The proof of (\ref{jed-b}) reduces to the proof of the fact that 
$Y^{1,p}_{t^\beta, c}\subseteq Y^{1,p}_{t^\beta ,0}$, which follows by the same convolution argument as just before.

\smallskip
\noindent
{\bf (\ref{jed-c}):}
To prove (\ref{jed-c}), we first show that $\ypb\subseteq 
\overline{Y^{1,p}_{t^\beta,c}}^{\,(\|.\|_{Y^{1,p}_{t^\beta}})}.$
For this, let $u\in \ypb$ and consider
\begin{equation}
\label{un}
u_n(t):= u(t)\phi_n(t)
\end{equation}
where
\begin{eqnarray*}
\phi_n(t):=\left\{
\begin{array}{lll}
0 &{\rm for}& t\le \frac{1}{2n} \vspace*{0.1cm}\\
\vspace*{0.1cm}
2nt-1 &{\rm for}& t\in (\frac{1}{2n},\frac{1}{n}] \vspace*{0.1cm}\\
1 &{\rm for}& t\in (\frac{1}{n},n] \vspace*{0.1cm}\\
-\frac{1}{n}t+2 &{\rm for}& t\in (n,2n] \vspace*{0.1cm}\\
0 &{\rm for}& t> 2n. 
\end{array}
\right.
\end{eqnarray*}
Clearly, $u_n(t)\in Y^{1,p}_{t^\beta,c}$. We will show that $u_n\to u$ in $Y^{1,p}_{t^\beta}$ as $n\to\infty$.
This follows from the following computations:
\begin{eqnarray}
 \int_{\rp} |u_n(t) -u(t)|^p t^\beta\left(1+\frac{1}{t^p}\right) dt &=& \int_{\rp} |u(t)(\phi_n(t) -1)|^p t^\beta\left(1+\frac{1}{t^p}\right) dt \nonumber\\
&\le &   \int_{ \{ (0,{\frac{1}{n}}) \cup (n,\infty)\} } |u(t)|^p t^\beta\left(1+\frac{1}{t^p}\right) dt  \stackrel{n\to\infty}{\rightarrow}0, \label{dwaosiem} \\
\int_{\rp} |\left( u_n(t) -u(t) \right)^{'}|^p t^\beta dt &=& \int_{\rp} |\left( u(t)(\phi_n(t) -1)\right)^{'}|^p t^\beta dt
\precsim  \int_0^{\frac{1}{n}} |u^{'}(t)|^p t^\beta dt \nonumber\\ +
\int_{n}^\infty |u^{'}(t)|^p t^\beta dt
&+& (2n)^p\int_{\frac{1}{2n}}^{\frac{1}{n}} |u(t)|^p t^\beta dt 
+\frac{1}{n^p}\int_n^{2n} |u(t)|^p t^\beta dt \nonumber \\
&=:& I_n+II_n+III_n + IV_n.  \label{granice1}
\end{eqnarray}
Obviously, $I_n,II_n,IV_n\to 0$ as $n\to \infty$. The same holds for $III_n$ because
\begin{eqnarray*}
III_n\sim \int_{\frac{1}{2n}}^{\frac{1}{n}} |u(t)|^p t^{\beta -p} dt
\leq \int_{\frac{1}{2n}}^{\frac{1}{n}} |u(t)|^p t^{\beta}\left(1+\frac{1}{t^p}\right) \stackrel{n\to\infty}{\rightarrow} 0.
\end{eqnarray*}

\smallskip
\noindent 
 It follows  that $$\ypb\subseteq 
\overline{Y^{1,p}_{t^\beta,c}}^{\,(\|.\|_{Y^{1,p}_{t^\beta}})}\stackrel{(\ref{jed-b})}{=} Y^{1,p}_{t^\beta,0}\subseteq  Y^{1,p}_{t^\beta} \ {\rm \text{and consequently,}}\ Y^{1,p}_{t^\beta,0}=Y^{1,p}_{t^\beta}.$$
Moreover, $$\ypb\subseteq 
\overline{Y^{1,p}_{t^\beta,c}}^{\,(\|.\|_{Y^{1,p}_{t^\beta}})}\subseteq 
\overline{W^{1,p}_{c}}^{\,(\|.\|_{Y^{1,p}_{t^\beta}})}\stackrel{(a)}{\subseteq} 
\overline{W^{1,p}_{c}}^{\,(\|.\|_{W^{1,p}(\rp , t^{\beta})})} \stackrel{ (\ref{jed-a})}{=} W_0^{1,p}(\rp , t^{\beta}),$$ 
where the embedding (a) holds because the norm  $\|\cdot\|_{W^{1,p}(\rp, t^{\beta})}$ is stronger than the norm  $\|\cdot\|_{\ypb}$.\hfill$\Box$

\section{Main results: analytic characterization of $W_0^{1,p}(\rp, t^\beta)$}
\subsection{Auxiliary sets and trace operators}
Let us define the following sets:
\begin{eqnarray*}
{\mathcal A}^{0}_{p,\,t^\beta}&:=&\{ u\in W^{1,p}(\rp , t^\beta) :  \exists_{t_n\searrow 0}: \lim_{n\to \infty} u(t_n)=0 \};\\
{\mathcal B}^{0}_{p,\,t^\beta}&:=&  \{ u\in W^{1,p}(\rp , t^\beta) : \lim_{t\to 0} u(t) =0\} ;\\
{\mathcal C}^{0}_{p,\,t^\beta}&:=&  \{ u\in W^{1,p}(\rp , t^\beta) :  \lim_{t\to 0} u(t)t^{(\beta -p+1)/p} =0\}; \\
{\mathcal A}^{\infty}_{p,\,t^\beta}&:=&\{ u\in W^{1,p}(\rp , t^\beta) :  \exists_{t_n \nearrow \infty}: \lim_{n\to \infty} u(t_n)=0 \};\\ 
{\mathcal B}^{\infty}_{p,\,t^\beta}&:=&  \{ u\in W^{1,p}(\rp , t^\beta) : \lim_{t\to \infty} u(t) =0\};\\
{\mathcal C}^{\infty}_{p,\,t^\beta}&:=&  \{ u\in W^{1,p}(\rp , t^\beta) : \lim_{t\to \infty} u(t)t^{(\beta -p+1)/p} =0\}.
\end{eqnarray*}
We will  consider the following two trace-endpoints operators, which are well defined on some subspaces $X, Z\subseteq W^{1,p}(\rp , t^\beta)$, possibly depending on $p$ and $\beta$:
\begin{equation}\label{tr}
Tr^0: X\rightarrow \BR,\ \hbox{\rm defined as}\ \   u\mapsto \lim_{t\to 0}u(t)
\end{equation}
and 
\begin{equation}
\label{tr2}
Tr^{\infty}: Z \rightarrow \BR , \ \hbox{\rm  defined as }\ u\mapsto \lim\limits_{t\rightarrow\infty}u(t).
\end{equation} 
In our further analysis we will indicate on  situations when  $X, Z$ are the whole spaces  $W^{1,p}(\rp , t^\beta)$. That observations will be useful to derive density results.

\subsection{Presentation of main results}

Our main goal is to prove the following three statements, which give  analytic characterization of $W_0^{1,p}(\rp, t^\beta)$ for all possible range of $\beta$.

\begin{theo}[characterization of $W_0^{1,p}(\rp, t^\beta)$, $\beta <p-1$]\label{pierwsze}
Let $1<p<\infty$ and $\beta <p-1$. Then we have:  
\begin{description} 
\item[i)] $W_0^{1,p}(\rp, t^\beta)= {\mathcal A}^{0}_{p,\,t^\beta}={\mathcal B}^{0}_{p,\,t^\beta}= {\mathcal C}^{0}_{p,\,t^\beta}  = 
\ypb$.  
\item[ii)]
The mapping $Tr^{0}$ as in \eqref{tr}
is well defined as the functional on $W^{1,p}(\rp , t^\beta)$. \\
Moreover,
\begin{equation}\label{trdalej}
\lim_{t\to 0} t^\frac{\beta- (p-1)}{p}|u(t)- Tr^0u| =0,\ \ \hbox{\rm for any}\  u\in W^{1,p}(\rp , t^\beta).
\end{equation}
\item[iii)] When $\beta\le -1$, then
$W^{1,p}(\rp , t^\beta)= W^{1,p}_0(\rp , t^\beta)={\mathcal B}^{0}_{p,\,t^\beta}  .$ In particular the trace operator $Tr^0$, defined by \eqref{tr}, annihilates the whole space $W^{1,p}(\rp , t^\beta)$
and 
\begin{equation}\label{trdalejjj}
W^{1,p}(\rp , t^\beta)= W_0^{1,p}(\rp, t^\beta)= {\mathcal A}^{0}_{p,\,t^\beta}={\mathcal B}^{0}_{p,\,t^\beta}= {\mathcal C}^{0}_{p,\,t^\beta} = 
\ypb .
\end{equation}
\item[iv)] When $\beta\in (-1,p-1)$, then the trace operator $Tr^0: W^{1,p}(\rp, t^\beta)\rightarrow \BR$ is an epimorphism. Moreover,
\begin{equation}\label{wtorek}
W^{1,p}(\rp , t^\beta)\varsupsetneq W_0^{1,p}(\rp , t^\beta)  
= {\mathcal A}^{0}_{p,\,t^\beta}={\mathcal B}^{0}_{p,\,t^\beta}= {\mathcal C}^{0}_{p,\,t^\beta}
=\ypb .
\end{equation}

\end{description}
\end{theo}

\begin{theo}[characterization of $W_0^{1,p}(\rp, t^\beta)$, $\beta > p-1$]
\label{trzecie}
Let $1<p<\infty$. When $\beta > p-1$, then we have:
 $$W^{1,p}(\rp, t^\beta)= W_0^{1,p}(\rp, t^\beta) = {\mathcal A}^{\infty}_{p,\,t^\beta}={\mathcal B}^{\infty}_{p,\,t^\beta}= {\mathcal C}^{\infty}_{p,\,t^\beta}=\ypb .$$
In particular the trace operator $Tr^\infty $, defined by formula \eqref{tr2}, annihilates the whole space  $W^{1,p}(\rp,t^\beta)$,  moreover
\begin{equation}
\label{trdalej2}
\lim_{t\to \infty } t^\frac{\beta- (p-1)}{p}|u(t)| =0,\ \ \hbox{\rm for any}\  u\in W^{1,p}(\rp , t^\beta).
\end{equation}
\end{theo}

\begin{theo}[characterization of $W_0^{1,p}(\rp, t^\beta)$,  $\beta =p-1$]\label{drugie}
Let $1<p<\infty$. Then we have:
\begin{description}
\item[i)] $W_0^{1,p}(\rp,  t^{p-1})= W^{1,p}(\rp, t^{p-1})$.
\item[ii)] $ Y^{1,p}_{t^{p-1}}\varsubsetneq W^{1,p}(\rp,  t^{p-1})$.
\item[iii)] The trace operator $Tr^0$ as in \eqref{tr} is not well defined on $W^{1,p}(\rp, t^{p-1})$. 
\item[iv)]
We have $W^{1,p}(\rp, t^{p-1}) = {\mathcal A}^{\infty}_{p,\,t^{p-1}}
 ={\mathcal B}^{\infty}_{p,\,t^{p-1}}$.
\end{description} 
  
\end{theo}

We have the following remark about the classical formulation of  Hardy inequality.
\begin{rem}\label{remhardy}\rm
 In the classical formulation of Hardy inequality (see Theorem \ref{hardyclas}), the sets of functions admitted to  Hardy inequality with right hand sides finite are
 $$\{ u\in W^{1,1}_{loc}(\rp): \lim_{t\to 0}u(t)=0, \int_{0}^\infty |u'(t)|^p\,t^{\beta}dt <\infty,  \}\ \rm{ when}\  \beta<p-1,$$  
$$\{  u\in W^{1,1}_{loc}(\rp): \lim_{t\to \infty}u(t)=0, \int_{0}^\infty |u'(t)|^p\,t^{\beta}dt <\infty, \}\ {\rm when}\  \beta>p-1.$$
We know from Theorems \ref{pierwsze} and \ref{trzecie}
that 
\begin{itemize}
\item
$W_0^{1,p}(\rp, t^\beta) = \{ u\in W^{1,p}(\rp, t^\beta): \lim_{t\to 0}u(t)=0\}$ when $\beta<p-1$,
\item
$W_0^{1,p}(\rp, t^\beta) = \{ u\in W^{1,p}(\rp, t^\beta): \lim_{t\to \infty}u(t)=0\}$ when $\beta>p-1$. 
\end{itemize}
Therefore the classical Hardy inequality, if restricted to $W^{1,p}(\rp, t^\beta)$, holds on the space $W_0^{1,p}(\rp, t^\beta)$. There is no larger subspace  $Z\subseteq W^{1,p}(\rp, t^\beta)$ admitted to the Hardy inequality with any finite constant. Indeed, for any $u\in W^{1,p}(\rp, t^\beta)\setminus W_0^{1,p}(\rp, t^\beta)$ we have either $\liminf_{t\to 0}|u(t)|=C>0$ or 
 $\liminf_{t\to \infty}|u(t)|=C>0$ by Theorems \ref{pierwsze} or \ref{trzecie}. But then $(|u|/t)^pt^\beta$ is not integrable either near zero or near infinity, so that left hand side in Hardy inequality is infinite. 
For recent improvement of Hardy inequality: 
$$\int_0^\infty {\rm max} \left\{ {\rm sup}_{0<s\le t}\frac{|u(s)|^p}{s^p}  ,
{\rm sup}_{t\le s <\infty}\frac{|u(s)|^p}{s^p}    \right\} dt \leq \left(  \frac{p}{p-1}\right)^p \int_0^\infty |u^{'}(t)|^p dt,$$
with the  admissibility condition $\liminf_{t\to 0} |u(t)|=0$, we refer to the  paper \cite{rulawe}. As the weight on the right hand side above is $1\equiv t^0$, according to Theorem \ref{pierwsze}, if we restrict 
such inequality to the subspace of $W^{1,p}(\rp)$, it must hold precisely on $W_0^{1,p}(\rp)$.
Hardy inequalities:
$$\left(\int_a^b |f(x)|^q u(x)dx\right)^\frac{1}{q}\leq C \left(\int_a^b |f^\prime(x)|^p v(x)dx\right)^\frac{1}{p},$$ where $-\infty\leq a<b\leq \infty,$ with
the analytical description of the admitted sets of functions given in the form of the vanishement conditions at the endpoints of the interval $(a, b)$ can also be found in \cite{kuf}. 
\end{rem}

In the proceeding section we will prove the above statements.

\subsection{The proofs}

\subsubsection*{Proof of Theorem \ref{pierwsze}.}

{\bf i) and ii):} The proof follows by steps.

\smallskip
\noindent
{\sc Step 1.} We show that ${\mathcal C}^{0}_{p,\,t^\beta} = {\mathcal B}^{0}_{p,\,t^\beta}= {\mathcal A}^{0}_{p,\,t^\beta}$. 
We obviously have
${\mathcal C}^{0}_{p,\,t^\beta}\subseteq {\mathcal B}^{0}_{p,\,t^\beta}\subseteq {\mathcal A}^{0}_{p,\,t^\beta}$. We will show that the converse inclusions hold.
Let $u\in {\mathcal A}^{0}_{p,\,t^\beta}$. Then there exists $t_n\to 0$ such that
$u(t_n)\to 0$. As for $0<t_n<t$
\begin{eqnarray}\label{tnt}
|u(t)-u(t_n)|\le \int_{t_n}^t |u^{'}(\tau )|\, d\tau \le 
\int_{0}^t |u^{'}(\tau )|\, d\tau ,
\end{eqnarray}
therefore,
\begin{eqnarray}\label{tnt1}
|u(t)|\le \int_{0}^t |u^{'}(\tau )|\, d\tau =  
 \int_{0}^t |u^{'}(\tau )|\tau^{\beta/p}\tau^{-\beta/p}\, d\tau
\\ 
\le 
\left( \int_0^t |u^{'}(\tau )|^p\tau^{\beta}\, d\tau\right)^{\frac{1}{p}} \left (\int_0^t \tau^{-\beta/(p-1)}\, d\tau \right)^{1-\frac{1}{p}}
= a(t)
t^{\frac{-\beta + p-1}{p}}\nonumber
\end{eqnarray}
where $a(t):=\left(\frac{p-1}{-\beta+p-1}\right)^\frac{p-1}{p}\left( \int_0^t |u^{'}(\tau )|^p\tau^{\beta}\, d\tau\right)^{\frac{1}{p}} \to 0$, as $t\to 0$. This implies
$$
t^{\frac{\beta - (p-1)}{p}}|u(t)| \stackrel{t\to 0}{\rightarrow} 0.
$$
Consequently, ${\mathcal A}^{0}_{p,\,t^\beta}\subseteq {\mathcal C}^{0}_{p,\,t^\beta}$, which ends the proof of Step 1.

\smallskip
\noindent
{\sc Step 2.} We prove {\bf ii)}.\\
Let us prove first that the operator $Tr^{0}$ as in (\ref{tr}) is well defined on $W^{1,p}(\rp , t^{\beta})$.

For this, let $u\in W^{1,p}(\rp , t^{\beta})$ and let $t_n\to 0$, as $n\to\infty$. 
By the same arguments as in (\ref{tnt}) and (\ref{tnt1}), we conclude that  
$$
|u(t) - u(t_n)|\precsim t^{-\frac{\beta-(p-1)}{p}}\left( \int_0^t |u^{'}(\tau)|^p \tau^\beta d\tau\right)^{\frac{1}{p}},\ {\rm for}\ t_n\le t.
$$
Therefore, the sequence  $\{ u(t_n)\}$ is bounded. By the Bolzano-Weierstrass theorem it possesses a convergent subsequence which we denote by the same expression.
We can assume that $u(t_n)\to \bar{u}\in \BR$ as $n\to\infty$.  
After passing to the limit in the above expression first with $n$, then with $t$, we get
\begin{eqnarray}\label{baruoszac}
|u(t) - \bar{u}|\precsim t^{-\frac{\beta-(p-1)}{p}}\left( \int_0^t |u^{'}(\tau)|^p\tau^\beta d\tau \right)^{\frac{1}{p}},\ \hbox{\rm  where}\ \bar{u}=\lim_{t_n\to 0} u(t_n), \\
\lim_{t\to 0} |u(t)-\bar{u}|=0\text{ and } \lim\limits_{t\to 0} t^\frac{\beta- (p-1)}{p}|u(t)-\bar{u}|=0.\nonumber
\end{eqnarray}

In particular, $u\mapsto Tr^0(u)=: u(0)$ is the well defined mapping on $W^{1,p}(\rp , t^{\beta})$. Such mapping is linear, so its continuity is equivalent to its boundedness.
In order to show the boundedness, we use H\"older inequality in the estimates below
\begin{eqnarray*}
|u(0)|= \int_0^1 |u(0)|d s \leq \int_0^1 |u(s)- u(0)|d s  + \int_0^1 |u(t)| t^{\frac{\beta}{p}}t^{-\frac{\beta}{p}} dt \\
\stackrel{(\ref{baruoszac})}\precsim 
 \int_0^1 s^{-\frac{\beta - (p-1)}{p}}\left( \int_0^1 |u^{'}(t)|^p
 t^\beta dt \right)^{\frac{1}{p}}d s  + 
\left( \int_0^1  |u(t)|^p t^\beta d t \right)^{\frac{1}{p}}\left(  \int_0^1 t^{-\beta/(p-1)} dt  \right)^{1-\frac{1}{p}}\\
 \precsim \left( \int_0^1 |u^{'}(t)|^p t^\beta dt\right)^{\frac{1}{p}} +
\left(  \int_0^1 |u(t)|^p t^\beta dt\right)^{\frac{1}{p}}\le \| u\|_{W^{1,p}(\rp , t^\beta)}.
 \end{eqnarray*}
This implies the continuity of the operator $Tr^0$ and allows to conclude that $Tr^0$ is the well defined functional on $W^{1,p}(\rp , t^\beta)$.
 Property (\ref{trdalej}) follows directly from (\ref{baruoszac}). 
 This ends the proof of Step 2.

\smallskip
\noindent
{\sc Step 3.} We show that ${\mathcal C}^{0}_{p,\,t^\beta}= W_0^{1,p}(\rp , t^\beta)$.\\
It follows from Step 2 that ${\mathcal B}^{0}_{p,\,t^\beta}$ can also be defined as $(Tr^0)^{-1}(0)$, so 
 ${\mathcal B}^{0}_{p,\,t^\beta} (={\mathcal C}^{0}_{p,\,t^\beta} = {\mathcal A}^{0}_{p,\,t^\beta})$ is a closed subspace in $W^{1,p}(\rp , t^\beta)$. Let $u\in {\mathcal B}^{0}_{p,\,t^\beta}$.
 
\smallskip
\noindent
By similar arguments as in the proof of formula (\ref{jed-c}) in Lemma 
\ref{gestosci},  defining $u_n$ ($n\in\BN$) by the  formula  (\ref{un}), we will show that ${\mathcal B}^{0}_{p,\,t^\beta}\subseteq \overline{W^{1,p}_{c}(\rp , t^\beta )}^{\,(\|.\|_{W^{1,p}(\rp , t^{\beta})})}.$
Indeed, the estimates as in (\ref{dwaosiem}) with the weight $t^\beta$ in place of $t^\beta (1+\frac{1}{t^p})$ show that $u_n\to u$ in $L^p(\rp , t^\beta )$. We repeat the estimates (\ref{granice1}) and then we only have to show that $III_n\to 0$ as $n\to\infty$. This is the consequence of the following estimates:
\begin{eqnarray}
III_n\sim \int_{\frac{1}{2n}}^{\frac{1}{n}} |u(t)|^p t^{\beta -p} dt
= \int_{\frac{1}{2n}}^{\frac{1}{n}} \left( 
|u(t)|^pt^{\beta - (p-1)} \right) t^{-1} dt \label{gwiazdka}\\
 \leq  {\rm sup}\{ 
|u(t)|^pt^{\beta - (p-1)}: t\le \frac{1}{n}\}\ln 2 \stackrel{n\rightarrow\infty}{\rightarrow}0.\nonumber
\end{eqnarray}
Last convergence holds because $u\in {\mathcal C}^{0}_{p,\,t^\beta} (={\mathcal B}^{0}_{p,\,t^\beta} )$
Having shown that $$
 C_0^\infty(\rp)\subseteq {\mathcal B}^{0}_{p,\,t^\beta}\subseteq 
\overline{W^{1,p}_{c}(\rp , t^\beta)}^{\,(\|.\|_{W^{1,p}(\rp , t^{\beta})})}\stackrel{Lemma\; \ref{gestosci}}{=} W_0^{1,p}(\rp , t^{\beta}),$$ and remembering ${\mathcal B}^{0}_{p,\,t^\beta}$ is the closed subset in $ W^{1,p}(\rp , t^{\beta})$, we deduce that
${\mathcal C}^{0}_{p,\,t^\beta}=W_0^{1,p}(\rp , t^{\beta})$, when passing to the respective closures. 
This finishes the proof of Step 3.

\smallskip
\noindent
{\sc Step 4.} We show that ${\mathcal B}^{0}_{p,\,t^\beta}=\ypb$.
\\
At first we observe that $\ypb\subseteq {\mathcal A}^{0}_{p,\,t^\beta}$.
Indeed, otherwise we would have $|u(t)|\ge C>0$ for the sufficiently small $t$, which implies the impossible inequality:
\[
\infty > \int_0^\varepsilon |u(t)|^pt^\beta (1+\frac{1}{t^p})dt \ge \int_0^\varepsilon |C|^p t^{\beta -p} dt =\infty . 
\]
Therefore $\ypb\subseteq W_0^{1,p}(\rp, t^\beta)$ by Steps 1 and 3. 
 Hardy inequality \eqref{clasha} implies that  $W_0^{1,p}(\rp ,t^\beta) (={\mathcal B}^{0}_{p,\,t^\beta})\subseteq \ypb$. This ends the proof of 
 Step 4 and part i). 

\smallskip
\noindent 
{\bf  iii):} When $\beta \le -1$ and $u\in L^p(\rp , t^\beta)$, we necessarily have $\liminf_{t\to 0} |u(t)|=0$, that is $u\in {\mathcal A}^{0}_{p,\,t^\beta}$. 
Indeed, otherwise we would have $|u(t)|\ge C>0$ for the sufficiently small $t$. Thus implies 
\begin{equation}
\label{zaba}
\infty> \int_0^\varepsilon |u(t)|^pt^\beta dt \ge \int_0^\varepsilon |C|^p t^\beta dt =\infty \text{ with some }\varepsilon>0.
\end{equation}

Thus, the condition $u\in W^{1,p}(\rp , t^\beta)$, together with
the already proven part i), gives that $u\in W_0^{1,p}(\rp ,t^\beta)={\mathcal B}^{0}_{p,\,t^\beta}$ and consequently $Tr^0u=0$.
Therefore, we get \eqref{trdalejjj}.

\smallskip
\noindent 
{\bf  iv):} 
Let $\phi \in Lip (\rp)$ be defined as 
\begin{eqnarray*}
\phi (t)=\left\{
\begin{array}{lll}
1 &{\rm  for} & t\in (0,1)\\
 -t +2& {\rm  for} & t\in [1,2)\\
0& {\rm for} & t\ge 2.
\end{array}
\right. 
\end{eqnarray*}
Obviously $\phi \in W^{1,p}(\rp , t^\beta)$ 
 and the extension operator 
 $$
\BR\ni s \mapsto Ext^0 (s):= s\phi (\cdot )\in W^{1,p}(\rp ,t^\beta)
$$
 is the right inverse to the operator of trace $Tr^0$. In particular the elements of  $W^{1,p}(\rp ,t^\beta)$ might have the nonzero limits at zero, so $W_0^{1,p}(\rp ,t^\beta)={\mathcal B}^{0}_{p,\,t^\beta}$ is the essential subspace of $W^{1,p}(\rp ,t^\beta)$. 
The identities (\ref{wtorek})
follow now from i).
 
This  finishes the proof of the  statement. 
 \hfill$\Box$

\subsubsection*{Proof of Theorem \ref{trzecie}.} 

 The equality of sets ${\mathcal A}^{\infty}_{p,\,t^\beta}={\mathcal B}^{\infty}_{p,\,t^\beta}= {\mathcal C}^{\infty}_{p,\,t^\beta}$ can be deduced by almost the same arguments as in the proof of subcase  i) Step 1 in Theorem \ref{pierwsze}. The main difference is that in the initial step we choose $T_n\to\infty$ instead of $t_n\to 0$ and we deal with integrals
$\int_t^{T_n}  |u^{'}(\tau )|\, d\tau$ and  
$\int_t^\infty |u^{'}(\tau )|\, d\tau$ $(t>0,t<T_n<\infty)$ instead of $\int_{t_n}^t |u^{'}(\tau )|\, d\tau$ and $\int_0^t |u^{'}(\tau )|\, d\tau$, respectively. 

Next, we observe that when $u\in L^p(\rp , t^\beta)$, then we necessarily have 
$$\liminf_{t\to \infty} |u(t)|=0,$$ 
as in the other case we would have:
$
\infty >\int_\varepsilon^\infty |u(t)|^pt^\beta dt =\infty$  where $\varepsilon >0$, 
by similar arguments as in \eqref{zaba}. 
Therefore, 
$$W^{1,p}(\rp, t^\beta)=  {\mathcal A}^{\infty}_{p,\,t^\beta}={\mathcal B}^{\infty}_{p,\,t^\beta}= {\mathcal C}^{\infty}_{p,\,t^\beta}. $$ 

\smallskip
The fact that $W^{1,p}(\rp, t^\beta)= W^{1,p}_0(\rp, t^\beta)$ follows from the similar computations as in \eqref{dwaosiem} (with the weight $t^\beta$) and \eqref{granice1}. When estimating $III_n$, we can use the fact that any $u\in W^{1,p}(\rp, t^\beta)$ converges to zero at $\infty$, so by the classical Hardy inequality  \eqref{clasha}, we have $\int_{\rp} \left(\frac{|u(t)|}{t}\right)^p t^\beta dt <\infty$. Consequently,
$$
III_n\sim \int_{\frac{1}{2n}}^{\frac{1}{n}} \left(\frac{|u(t)|}{t}\right)^pt^\beta dt \stackrel{n\to\infty}{\rightarrow} 0.
$$

Note also that obviously $\ypb \subseteq W^{1,p}(\rp , t^\beta)$, while  Hardy inequality \eqref{clasha}, together with the fact that 
$W^{1,p}(\rp, t^\beta)={\mathcal B}^{\infty}_{p,\,t^\beta}$
 implies that $\ypb \supseteq W^{1,p}(\rp , t^\beta)$.

\smallskip
\noindent {\bf (\ref{trdalej2}):} The property \eqref{trdalej2} follows from modification of \eqref{tnt1}
(now we choose $T_n\to\infty$ and the integral over the integral $(t,\infty)$ instead of $(0,t)$), which reads as:
$$
|u(t)|\le t^{-\frac{\beta - (p-1)}{p}}\bar{a}(t)
$$
where
 $$\bar{a}(t):=\left(\frac{p-1}{\beta-(p-1)}\right)^\frac{p-1}{p}\left( \int_t^\infty |u^{'}(\tau )|^p\tau^{\beta}\, d\tau\right)^{\frac{1}{p}} \to 0,\ {\rm as}\  t\to \infty,$$ 
which finishes the proof of the statement. \hfill$\Box$

\subsubsection*{Proof of Theorem \ref{drugie} }

{\bf i):} 
The proof follows by steps.

\smallskip
\noindent
{\sc Step 1.} We show  that ${\mathcal B}^{0}_{p,\,t^{p-1}} \subseteq W_0^{1,p}(\rp, t^{p-1})$. 

For this, let  $u\in {\mathcal B}^{0}_{p,\,t^{p-1}}$ and  
 consider the sequence $\{ u_n\}_{n\in\mathbb{N}}$ as in (\ref{un}).
Similar estimates as in \eqref{dwaosiem} (with the weight $t^{p-1}$)  and
 (\ref{gwiazdka}) show that $u\in W_0^{1,p}(\rp, t^{p-1})$.
 
\smallskip
\noindent 
{\sc Step 2.} We show  that the  set $L\subseteq  W_0^{1,p}(\rp, t^{p-1})$ where
$$L:= \{ u\in W^{1,p}(\rp , t^{p-1}) : \ \hbox{\rm $u$ is bounded near zero}\}\ \  (\supseteq {\mathcal B}^{0}_{p,\,t^{p-1}} ) .$$

 For this, let us we consider $u\in L$. By the same computations as in Step 1, we deduce that the sequence $\{ u_n\}$ as in (\ref{un}), is bounded in $W^{1,p}(\rp, t^{p-1})$. This is  because by (\ref{gwiazdka}), we have
$$III_{n}\precsim \| u\|_{L^\infty (0,1)}\ln 2.$$ 
The argument as in \eqref{dwaosiem} with the weight $t^{p-1}$
 gives that $u_n\to u$ strongly in $L^p(\rp, t^{p-1})$. We will prove that 
\begin{equation}\label{pochodne}
u_n^{'}\stackrel{n\to \infty}{\rightharpoonup} u^{'}\ \hbox{\rm  (weakly) in}\   L^p(\rp, t^{p-1}).
\end{equation}

To verify (\ref{pochodne}) at first we note that the space $L^{p/(p-1)}(\rp, t^{-1})$   is  dual to $L^{p}(\rp, t^{p-1})$, where the duality is expressed by the formula
\begin{eqnarray}
L^p(\rp, t^{p-1})\times L^{p/p-1}(\rp, t^{-1}) \ni (w,v) & \mapsto 
\int_{\rp} w(t) v(t) dt= \int_{\rp} (w\cdot t^{1-\frac{1}{p}})\cdot(v\cdot t^{\frac{1}{p}-1}) dt. \nonumber
\end{eqnarray}
Therefore it suffices to test the weak convergence \eqref{pochodne} on any dense subset in \sloppy $L^{p/(p-1)}(\rp, t^{-1})$, such as the set of smooth compactly supported functions. Hence, it suffices to verify that for any $v\in C_0^\infty(\BR_+)$ supported in $[a,b]$, where $0<a<b<\infty$, we have
$$
\int_{\rp} (u_n^{'}(t)-u^{'}(t))v(t) dt \stackrel{n\to\infty}{\rightarrow} 0. 
$$
This follows from the following simple argument:
\begin{eqnarray*}
\int_{\rp} (u_n^{'}(t)-u^{'}(t))v(t) dt = \int_{[a,b]} (u_n^{'}(t)-u^{'}(t))v(t) dt\\ = \int_{[a,b]} u^{'}(t)(\phi_n(t)-1)v(t) dt + 
\int_{[a,b]} u(t)\phi_n^{'}(t)v(t) dt =: I_n + II_n.
\end{eqnarray*}
As $\phi_n \equiv 1$ on sets $[\frac{1}{n}, n]$, which exhaust $[a,b]$ for the sufficiently  large $n$, therefore for such $n$ we have $I_n = II_n =0$.

Mazur's theorem (see, e.g. \cite{ek}, p. 6) says that there exists the sequence $\{ v_n\}$ of convex combinations of the $u_n$'s such that $v_n^{'}\rightarrow u^{'}$ in $L^p(\rp, t^{p-1})$.
Simple verification shows that also $v_n\to u$ strongly in $L^p(\rp,t^{p-1})$. In particular $C_0^\infty(\rp)\ni v_n\to u$ in $W^{1,p}(\rp, t^{p-1})$, as $n\to\infty$ and consequently, $u\in W_0^{1,p}(\rp, t^{p-1})$.

\smallskip
\noindent 
{\sc Step 3.} We show that  $W^{1,p}(\rp, t^{p-1})=  W_0^{1,p}(\rp, t^{p-1})$, which will finish the proof of part i).

For this, we note that $$W_0^{1,p}(\rp , t^{p-1})\stackrel{Step\  2}{\supset }L\supset L^\infty (\rp) \cap W^{1,p}(\rp , t^{p-1}).$$
Moreover, $L^\infty (\rp) \cap W^{1,p}(\rp , t^{p-1})$ is a dense subset in $W^{1,p}(\rp , t^{p-1})$. This is  because the truncation of any $u\in W^{1,p}(\rp , t^{p-1})$:
\begin{eqnarray*}
T_\lambda (u)(t):= 
\left\{
\begin{array}{lll}
u(x) & {\rm if}& |u(x)|\le \lambda\\
-\lambda & {\rm if} & u(x)\le -\lambda\\
\lambda & {\rm if} & u(x)\ge \lambda\\
\end{array}
\right.
\end{eqnarray*}
converge to $u$, as $\lambda \to\infty$, in  $W^{1,p}(\rp , t^{p-1})$ because
\begin{eqnarray*}
\int_{\BR_+} |u-T_\lambda(u)|^p t^{p-1} dt & = & \int_{u>\lambda} |u-\lambda |^pt^{p-1}dt + \int_{u< -\lambda} |u+\lambda |^p t^{p-1} dt \\
& = & \int_{|u|>\lambda} ||u|-\lambda |^pt^{p-1}dt \leq \int_{|u|>\lambda} |u|^p t^{p-1}dt \overset{\lambda\rightarrow\infty}{\rightarrow}0
\end{eqnarray*}
and similarly
\begin{eqnarray*}
\int_{\BR_+} |u^\prime -(T_\lambda(u))^\prime |^p t^{p-1} dt & = & \int_{|u|>\lambda} |u^\prime |^p t^{p-1} dt \overset{\lambda\rightarrow\infty}{\rightarrow}0.
\end{eqnarray*}

Therefore, $W^{1,p}(\rp , t^{p-1})\subseteq  W_0^{1,p}(\rp , t^{p-1})$, which ends the proof of Step 3.

\smallskip
\noindent
{\bf ii):}  As $u(t) = \ln (-\ln t )\chi_{(0,\frac{1}{e})}\not\in Y^{1,p}_{t^{p-1}}$, therefore 
$Y^{1,p}_{t^{p-1}} \varsubsetneq W^{1,p}(\rp, t^{p-1})$. 

\smallskip
\noindent
{\bf iii):} The function 
$u(t) = \ln (-\ln t )\chi_{(0,\frac{1}{e})}$ 
belongs to $W^{1,p}(\rp, t^{p-1})$ but converges to $\infty$ at zero. 
Therefore, the trace operator as in (\ref{tr}) cannot be well defined on the space $W^{1,p}(\rp, t^{p-1})$.

\smallskip
\noindent
{\bf iv):} The fact that $W^{1,p}(\rp , t^{p-1})\subseteq {\mathcal A}^{\infty}_{p,\,t^{p-1}}$ 
follows from integrability of $|u|^pt^{p-1}$ near infinity and the argument as in \eqref{zaba}. We are left with the proof that 
$${\mathcal A}^{\infty}_{p,\,t^{p-1}}= {\mathcal B}^{\infty}_{p,\,t^{p-1}}.$$
Suppose on the contrary, that the above is not true, that is there exists
the function $u\in W^{1,p}(\rp, t^{p-1})$, which is oscillating near infinity, that is $\limsup_{t\to\infty}|u(t)|=c>0=\liminf_{t\to\infty} |u(t)|$. Using the Darboux property and changing the function $u$ to $Au$ with the suitable chosen constant $A$, we could construct the function 
$v\in W^{1,p}(\rp, t^{p-1})$, and sequences $t_n\nearrow \infty$, $T_n \nearrow \infty$ such that
$$ 0<t_{n-1}< T_n < t_n \ {\rm and}\ v(t_n)=\frac{1}{2},\ v(T_n)=1,\ \frac{1}{2}< v(t)<1 \ {\rm on}\ (t_{n-1},T_n).$$
This would imply:
\begin{eqnarray}
\infty &>&\int_0^\infty |v(t)|^pt^{p-1}dt \ge  \sum_{n=1}^\infty \int_{t_{n-1}}^{T_n}\left( \frac{1}{2}\right)^p t^{p-1}dt  
 =   \frac{1}{p2^p}   \sum_{n=1}^\infty
\left( T_n^p-t_n^p\right)\nonumber\\ &=& \frac{1}{p2^p}   \sum_{n=1}^\infty
T_n^p\left( 1-\left(\frac{t_{n-1}}{T_n}\right)^p\right)
=: \frac{1}{p2^p}   \sum_{n=1}^\infty
T_n^p\left( 1-(x_n)^p\right). \label{marzec1}
\end{eqnarray}
At the same time 
\begin{eqnarray*} 
\frac{1}{2}&= &|v(T_n)-v(t_n)| \le \int_{t_{n-1}}^{T_n} |v^{'}(t)|dt =
\int_{t_{n-1}}^{T_n} |v^{'}(t)|t^{\frac{p-1}{p}}t^{-\frac{p-1}{p}} dt
 \nonumber\\
 &\le &
\left( \int_{t_{n-1}}^{T_n} |v^{'}(t)|^pt^{p-1}dt\right)^{\frac{1}{p}}
\left( \int_{t_{n-1}}^{T_n} t^{-1}dt\right)^{1-\frac{1}{p}}=: a_n b_n.
  \end{eqnarray*}
As $a_n\to 0$ when $n\to\infty$, therefore $b_n\to\infty$ when $n\to\infty$.
Consequently, $b_n = \ln (\frac{T_n}{t_{n-1}})\to \infty$ as $n\to\infty$ and 
so $\{ x_n\}$ in \eqref{marzec1} satisfies $x_n=\frac{t_{n-1}}{T_n}\to 0$ as $n\to \infty$, in particular $(1-x_n^p)>\frac{1}{2}$ for the sufficiently large $n$.\\
Further estimates in \eqref{marzec1} would lead to the impossible inequality, with some number  $N\in\BN$:
$$
\infty > \int_0^\infty |v(t)|^pt^{p-1}dt 
\succsim \sum_{n= N}^\infty T_n^p = \infty .
$$
We arrive at contradiction, which confirms that ${\mathcal A}^{\infty}_{p,\,t^{p-1}}= {\mathcal B}^{\infty}_{p,\,t^{p-1}}$ and ends the proof of the statement. 
\hfill$\Box$

\subsection{Applications}

We propose three example applications of our results.

\begin{rem}[application to BVP]\label{applbvp}\rm
Our results can be applied to the analysis of Boundary Value Problems in O.D.E's. For example, when dealing with the concrete ODE with the solution in $W_0^{1,p}(\rp, t^\beta)$, where $p>1$, one can formulate the boundary decay condition as:
 \begin{eqnarray*}
 \lim_{t\to 0} u(t)= 0, \ {\rm equivalently}\  \lim_{t\to 0} t^{(\beta-(p-1))/p} u(t),\ \ \hbox{\rm when}\ \beta <p-1,\\
  \lim_{t\to \infty} u(t)= 0, \ {\rm equivalently}\ \lim_{t\to \infty} t^{(\beta-(p-1))/p} u(t),\ \ \hbox{\rm when}\ \beta >p-1,
  \end{eqnarray*}
for  $u\in W_0^{1,p}(\rp, t^\beta)$. 
In particular, all the above limits are equivalent to the limits of $t^{s} u(t)$, for any $s\in [0, (\beta-(p-1))/p]$. 
We have shown in Theorems \ref{pierwsze} and \ref{trzecie} that such limits are well defined in the respective Sobolev spaces. Such limits can be also prescribed to the radial solutions to the Boundary Value Problems in PDE's defined on $\mathbb{R}^n$, where we expect radial solutions $v(x)=u(\|x\|)$ in the weighted Sobolev Spaces $W^{1,p}(\mathbb{R}^n, \| x\|^\alpha)$, where $\beta =\alpha +n-1$. 
\end{rem}

\begin{rem}[application to Hardy-type inequality]\label{appl-interpo2}\rm
In \cite{akikppcejm} K. Pietruska-Pa\l{}uba and second author have studied the inequalities:
\begin{equation}\label{zabak}
\int_{\rp} M(t^\alpha |u(t)|)t^\beta dt \le C_1 \int_{\rp} M(|u(t)|)t^\beta dt + C_2\int_{\rp} M(|u^{'}(t)|)t^\beta dt,
\end{equation}
holding on supersets $W$  of $C_0^\infty (\rp)$, where $M(\cdot)$ is the given Orlicz function, as the special cases of larger class of inequalities 
from \cite{akikppmia12}, which deal  with general weights $\omega, \rho$ in place of $t^\alpha, t^\beta$, respectively.
Let us focus on  \eqref{zabak} restricted to $M(s)=s^p, p>1$:
\begin{equation}\label{zabakp}
\int_{\rp} |t^\alpha |u(t)||^pt^\beta dt \le C_1 \int_{\rp} |u(t)|^pt^\beta dt + C_2\int_{\rp} |u^{'}(t)|^pt^\beta dt.
\end{equation}
 The following statement is the special case of Theorem 4.2 from 
\cite{akikppcejm}, where we put   $M(s)=s^p$. 
\begin{theo}
Suppose that $\alpha\in (-1,0)$. Then we have:
\begin{description}
\item[i)] Inequality  \eqref{zabakp} with $C_1=0$ cannot hold for all $u\in C_0^\infty(\rp)$ with any finite constant independent on $u$.
\item[ii)] There exist positive constants $C_1,C_2$ such that inequality
\eqref{zabakp} holds for all $u\in W(\alpha,\beta)=W$, where\\
a) $W=W^{1,p}(\rp, t^\beta)$ when $\beta >|\alpha|p-1$,
\\
b) $W=\{ u\in W^{1,p}(\rp, t^\beta) : \liminf_{t\to 0^+}|u(t)|t^s =0, \ s:={\frac{\beta+1-|\alpha | p}{p}  } \}$ when $\beta <|\alpha|p-1$.
\item[iii)] The sets $W(\alpha,\beta)$ from ii)  are maximal subsets of 
$W^{1,p}(\rp, t^\beta)$ for which \eqref{zabakp} holds.
\end{description}
\end{theo}
We can now show that when $\beta <|\alpha|p-1$ then 
$W=W_0^{1,p}(\rp, t^\beta)$ and in particular it is independent on $\alpha$.
Indeed, for any $u\in W=W(\alpha,\beta)$, we necessarily have 
$\liminf_{t\to 0}|u(t)|=0$, because $s<0$. Therefore 
$u\in \mathcal{A}^0_{p,t^\beta}$.
As $\beta <p-1$, therefore by Theorem \ref{pierwsze} we get $u\in W_0^{1,p}(\rp, t^\beta)$. This gives $W \subseteq W^{1,p}_0(\rp, t^\beta)$.
On the other hand, when $u\in W_0^{1,p}(\rp, t^\beta)$ then  by Theorem \ref{pierwsze} $u\in \mathcal{C}^0_{p,t^\beta}$ and consequently
\begin{eqnarray*}
\lim_{t\to 0} |u(t)|t^{\frac{\beta+1-p}{p}}=0 \Longrightarrow
\lim_{t\to 0} |u(t)|t^{\frac{\beta+1-|\alpha |p}{p}}=
 \lim_{t\to 0} \left(|u(t)|t^{\frac{\beta+1-p}{p}}  \right) t^{1-|\alpha|} =0. 
\end{eqnarray*}
Hence $u\in W$. We have shown that $W=W_0^{1,p}(\rp, t^\beta)$.
\end{rem}

\begin{rem}[application to complex interpolation theory]\label{appl-interpo}\rm
It is known by Stein - Weiss Theorem (see i.e. \cite{bl}, Section 5.5.3, p. 120) that having two weight functions $\omega_0,\omega_1$ defined on some domain $U\subseteq \mathbb{R}^n$, one has
$$
[ L^{p_0}(U,\omega_0, \mathbb{C}), L^{p_1}(U,\omega_1,\mathbb{C})]_\theta =L^{p_\theta }(U,\omega_\theta, \mathbb{C}),\ \ \hbox{\rm where}\ \theta \in (0,1), 
$$
where $[X,Y]_\theta$ denotes the complex interpolation pair, 
\begin{equation}
\label{dwaw}
\frac{1}{p_\theta} =\frac{1-\theta}{p_0} + \frac{\theta}{p_1},\;\;\omega_\theta^{\frac{1}{p_\theta}} = \omega_0^{\frac{1-\theta}{p_0}}\omega_1^{\frac{\theta}{p_1}}
\end{equation}
and $L^{p}(U,\omega_0, \mathbb{C})= L^{p}(U,\omega_0)\times iL^{p}(U,\omega_0)$  denotes complex valued functions $f_1+if_2$, with coordinate functions $f_1,f_2$  from
$L^{p}(U,\omega_0)$.

Let us use the analogous notation for complex-valued Sobolev functions 
$W^{s,p}(U,\omega,\mathbb{C})=  W^{s,p}(U,\omega )\times iW^{s,p}(U,\omega ) $. As proven by J\"orgen L\"ofstr\"om \cite{Lo}, under some special assumptions on weight functions, one has
$$
[ W^{s_0,p_0}(\BR^n,\omega_0,\mathbb{C}), W^{s_1,p_1}(\BR^n,\omega_1,\mathbb{C})]_\theta =W^{s_\theta,  p_\theta }(\BR^n,\omega_\theta,\mathbb{C}),
$$
where $\theta \in (0,1)$ 
and $s_\theta =(1-\theta)s_0+\theta s_1$, $p_\theta$ and $\omega_\theta$ are given by \eqref{dwaw} and weighted function Sobolev spaces are defined with the help of Fourier transform. The authors of \cite{ce} have raised  interesting questions, when one has
$$
[ W^{1,p}(U,\omega_0,\mathbb{C}), W^{1,p}(U,\omega_1,\mathbb{C})]_\theta =W^{1,  p }(U,\omega_\theta,\mathbb{C}),
$$
or
$$
[ W_0^{1,p}(U,\omega_0,\mathbb{C}), W_0^{1,p}(U,\omega_1,\mathbb{C})]_\theta =W_0^{1,  p }(U,\omega_\theta,\mathbb{C}),
$$
 where $\theta \in (0,1)$, 
with the same $\omega_\theta$ as in \eqref{dwaw} and $p\in [1,\infty)$?

It has been shown in \cite{ce} in Theorem 1.16 that under the assumption (1.6) on weights, which are satisfied for positive continuous functions, and local Lipschitzism on $r_\omega$ one has
\begin{eqnarray}
[W^{1,p}(U,\omega_0,\mathbb{C}), W^{1,p}(U,\omega_1, \mathbb{C})]_\theta  \subseteq W^{1,p}(U,\omega_\theta,  \mathbb{C}), \label{zero1}~~~~~~~~~~~~~~~~~~~~~~~~~\\
\mathcal{W}_0^p(U,\theta, r_w)\subseteq  [ W_0^{1,p}(U,\omega_0, \mathbb{C}), W_0^{1,p}(U,\omega_1, \mathbb{C})]_\theta ,~~~~~~~~~~~~~~~~~~~~~~~~\label{zero}\\
\mathcal{W}^p(U,\theta, r_w) :=  \{ \phi :U\rightarrow \mathbb{C} : \phi \in W^{1,p}(U, w_\theta,\mathbb{C}), \phi\cdot |\nabla \log(r_w)|\in L^p(U, w_\theta , \mathbb{C}) \} ,\nonumber\\
\|\phi\|_{\mathcal{W}^p(U,\,\theta,\, r_w)}: =\left(\|\phi\|_{W^{1,p}(U,\,\omega_\theta, \mathbb{C})}^p + \int_U |\phi(x)|^p|\nabla\log(r_\omega(x))|^p\omega_\theta(x)dx \right)^\frac{1}{p}, \ \    
 r_\omega : = \frac{\omega_0}{\omega_1}, \nonumber
\end{eqnarray} 
$\mathcal{W}_0^p(U,\theta, r_w)$ is the completion of Lipschitz compactly supported complex-valued functions in $\mathcal{W}^p(U,\theta, r_w)$.
As results from \eqref{zero1} and \eqref{zero}, we have:
\begin{eqnarray*}
\mathcal{W}_0^p(U,\theta, r_w)&\subseteq &  [ W_0^{1,p}(U,\omega_0, \mathbb{C}), W_0^{1,p}(U,\omega_1, \mathbb{C})]_\theta\stackrel{(a)}{\subseteq}  [W^{1,p}(U,\omega_0,\mathbb{C}), W^{1,p}(U,\omega_1, \mathbb{C})]_\theta  \\
&\subseteq & W^{1,p}(U,\omega_\theta,  \mathbb{C}),
\end{eqnarray*}
where the inclusion (a) follows from the very definition of interpolation pair, see Remark 1.17 in \cite{ce}.
Therefore if one has 
\begin{equation*}
\mathcal{W}_0^p(U,\theta, r_w) = W^{1,p}(U,\omega_\theta,  \mathbb{C}),
\end{equation*}
then 
\begin{equation*}
[ W_0^{1,p}(U,\omega_0, \mathbb{C}), W_0^{1,p}(U,\omega_1, \mathbb{C})]_\theta=  [W^{1,p}(U,\omega_0,\mathbb{C}), W^{1,p}(U,\omega_1, \mathbb{C})]_\theta  
= W^{1,p}(U,\omega_\theta,  \mathbb{C})
\end{equation*} 

Our main results cannot be concluded, as a special case, from Theorem 4 in \cite{Lo} because our weights $t^\beta$ do not satisfy the "polynomially regularity" condition, an essential and important assumption for the  weight function, considered among others in Theorem 4 of \cite{Lo}. 
Even if our functions could be extended to those on the whole $\BR$, it is not possible to provide further analysis based on paper \cite{Lo} because the polynomial components of the weight constrain the fulfillment of condition $(3.4)$ from the definition of polynomially regular weight function in \cite{Lo} on page 197, which requires in our case the finiteness of $\sup\limits_{t\in\BR_+}\vert\frac{(t^\beta)^\prime}{t^\beta}\vert=\sup\limits_{t\in\BR_+}\frac{\vert\beta\vert}{t}$. However, in our case this supremum is infinite.

\smallskip
\noindent 

Let us consider $U=\rp$, $w_0=t^{\beta_0}$, $\omega_1=t^{\beta_1}$. Then we have
$$
r_\omega(t)= t^{\beta_0-\beta_1}, \ \omega_\theta  = t^{\beta_0\theta + (1-\theta)\beta_1}=:t^{\beta_\theta},\ \ |\nabla \log(r_\omega)|\sim \frac{1}{t}
$$
and
\begin{eqnarray*}
\mathcal{W}^p(\rp ,\theta, t^{\beta_0-\beta_1}) 
=\\
\{ \phi : \rp\rightarrow\mathbb{C}: \phi \sim (\phi_1,\phi_2),   \ 
 \phi_i \in W^{1,p}(\rp , t^{\beta_{\theta}})\ {\rm for}\ i=1,2,\ 
 \phi/t \in L^p (\rp, t^{\beta_\theta})\}\\ = 
 Y^{1,p}_{t^{\beta_\theta}}\times i Y^{1,p}_{t^{\beta_\theta}} . 
\end{eqnarray*}
It follows from Theorem \ref{pierwsze} that 
\begin{eqnarray*}
W_0^{1,p}(\rp, t^{\beta_\theta},\mathbb{C}) = W^{1,p}(\rp, t^{\beta_\theta}, \mathbb{C}), 
 {\rm equivalently}\ W_0^{1,p}(\rp, t^{\beta_\theta}) = W^{1,p}(\rp, t^{\beta_\theta})\\ \Longleftrightarrow \beta_\theta \in (-\infty , -1]\cup [p-1,+\infty),  
 \end{eqnarray*}
while  from \eqref{jed-b} in Lemma \ref{gestosci} we know that
\begin{eqnarray*}
\mathcal{W}^p(\rp ,\theta, t^{\beta_0-\beta_1})= \mathcal{W}_0^p(\rp ,\theta, t^{\beta_0-\beta_1}) \\
\Longleftrightarrow 
Y^{1,p}_{t^{\beta_\theta}}= Y^{1,p}_{t^{\beta_\theta}, 0},\ \   \hbox{\rm which holds for any}\ \beta_\theta\in \mathbb{R}, p>1.
\end{eqnarray*}
Moreover, when only $p>1, \beta_\theta\neq p-1$, then by Theorem \ref{pierwsze}
$$ 
\mathcal{W}^p(\rp ,\theta, t^{\beta_0-\beta_1})= W_0^{1,p}(\rp, t^{\beta_\theta},\mathbb{C}), \ {\rm equivalently} \ Y^{1,p}_{t^{\beta_\theta}} = W_0^{1,p}(\rp, t^{\beta_\theta}).$$
Altogether leads to the following interpolation result, which seems to us new.

\begin{theo}\label{interpo-konkluzja-my}
When $\theta\in (0,1)$ $\beta_0,\beta_1\in\mathbb{R}$ are such that for
$\beta_\theta =(1-\theta) \beta_0 +\theta\beta_1\in (-\infty, -1]\cup (p-1,\infty)
$, 
then 
\begin{eqnarray*}
[ W_0^{1,p}(\rp, t^{\beta_0}, \mathbb{C}), W_0^{1,p}(\rp,t^{\beta_1}, \mathbb{C})]_\theta=  [W^{1,p}(\rp,t^{\beta_0},\mathbb{C}), W^{1,p}(\rp,t^{\beta_1}, \mathbb{C})]_\theta \nonumber \\
= W^{1,p}(\rp ,t^{\beta_\theta},  \mathbb{C})= Y^{1,p}_{t^{\beta_\theta}}\times i Y^{1,p}_{t^{\beta_\theta}}.
\end{eqnarray*} 

\end{theo}
The related statement  is Theorem 1.20 in \cite{ce},   which gives the condition for 
$$\mathcal{W}_0^p(U,\theta, r_w) \subseteq [ W_0^{1,p}(\rp, t^{\beta_0}, \mathbb{C}), W_0^{1,p}(\rp,t^{\beta_1}, \mathbb{C})]_\theta \subseteq W_0^{1,p}(\rp ,t^{\beta_\theta},  \mathbb{C})$$
 for all $\theta\in (0,1)$, in the abstract approach. For that inclusions one has to verify the following assumption: 
 
 {\it "whenever $\phi$ is an element of $W_0^{1,p}(U,\omega_0)\cap W_0^{1,p}(U,\omega_1)$ there exists a sequence $\{ \phi_n\}$ of functions in $Lip_c(U)$ which converges to $\phi$ in $W^{1,p}(U,\omega_0)$ or in $W^{1,p}(U,\omega_1)$ and is bounded in the other space."} 

It shows that the knowledge about density results and the techniques of approximation can be the crucial  tool for the analysis of interpolation spaces in the complex interpolation theory.
\end{rem}

\end{document}